\documentclass{amsart}
\usepackage{amssymb}
\newtheorem{Theo}{Theorem}

\newtheorem{Con}[Theo]{Conjecture}
\newtheorem{Lem}[Theo]{Lemma}

\begin{document}
\title{On triangular billiards}
\author{Jan-Christoph Schlage-Puchta}
\begin{abstract}
We prove a conjecture of Kenyon and Smillie concerning the nonexistence of
acute rational-angled triangles with the lattice property.
\end{abstract}

\maketitle

In a recent paper\cite{KS} on Billiards on rational-angled triangles, R.
Kenyon and J. Smillie proved the following theorem:

\begin{Theo}
Let $T$ be an acute non-isosceles rational angled triangle with angles
$\alpha$, $\beta$ and $\gamma$, which can be written as $p_1\pi/q$, $p_2\pi/q$
and $p_3\pi/q$ with $q\leq 10000$. Then $T$ is a polygon with the lattice
property if and only if
$(\alpha, \beta, \gamma)$ is one of the following:
\[
(\pi/4, \pi/3, 5\pi/12),\quad (\pi/5, \pi/3, 7\pi/15),\quad
(2\pi/9, \pi/3, 4\pi/9).
\]
\end{Theo}

They further showed, that the restricition on $q$ may be dropped, if the
following conjecture was true(see \cite{KS}, p. 94f):

\begin{Con}
Let $n, s, t$ be integers with $(n, s)=1$, $1\leq s, t<n$. Assume that for all
$p$ with $(p, n)=1$ we have
$\frac{n}{2}<ps\bmod n + pt\bmod n < \frac{3n}{2}$. Then
one of the following conditions hold true: $n\leq 78$, $s+t=n$, $s+2t=n$,
$2s+t=n$, or $n$ is even, and $|t-s|=\frac{n}{2}$.
\end{Con}

In this note we will prove this conjecture:

\begin{Theo}
Conjecture 2 is true.
\end{Theo}

Note that the classification of non-obtuse rational angled triangles with the
lattice-property is complete, since the cases of isosceles and right angled
triangles are completely solved in \cite{KS}, too.

By direct calculation, R. Kenyon and J. Smillie showed, that Theorem 3 is true
for $n\leq 10000$. We will use this fact at several steps in the proof.

The proof will depend on several facts concerning the distribution of relative
prime residue classes, collected in the next Lemma. We write $g(n)$ for the
Jacobsthal function, given by the maximal difference of consecutive integers
relatively prime to $n$, and $\omega(n)$ for the number of distinct prime
factors of $n$.

\begin{Lem}
\begin{enumerate}
\item We have $g(n)\leq 2^{\omega(n)}$. If $\omega(n)\leq 12$, we have
$g(n)\leq \omega(n)^2$.
\item Assume that $(a, d, n)=1$. Then in every interval $[x, x+g(n)]$ there is
some integer $\nu$, such that $(n, d\nu+a) = 1$.
\item For all $d>2$ there exists some $a$ with $(d, a)=1$ and
$\frac{d}{12}<a<\frac{5d}{12}$.
\item If $m$ is the product of the first $\omega(n)$ prime numbers, then
$g(n)\leq g(m)$.
\item We have $g(30)=6$, $g(210) = 10$, $g(2310) = 14$, $g(30030)= 22$,
$g(510510)=26$, $g(9699690) = 34$.
\end{enumerate}
\end{Lem}
{\em Proof:} The first statement was proven by Kanold\cite{Kan}. To prove the
second statement note first that it is trivial if $(d, n)=1$, for if
$dd'\equiv 1\pmod{n}$, then
the integers $dd'\nu+d'a$ are consecutive $\pmod{n}$, and none is coprime to
$n$, contradicting the definition of $g$. Now we may without loss assume that
$n$ is squarefree. If $(d, n) = e>1$, the integers $d\nu+a$ are coprime to $n$
if and only if they are coprime to $n/e$, thus using the case $(n, d)=1$ we get
that there is some $\nu\in[x, x+g(n/e)]$ such that $(d\nu+a, n)=1$. 
The third statement follows for $d>30$ from the first one, for $3\leq d\leq 30$
by direct inspection. The fourth statement was proven by Iwaniec\cite{Iwa}.
The fifth statement can be checked by direct computation.

Note that the fourth and fifth statement together greatly improve the
first one for $\omega(n)\leq 8$.

Note further that the asymptotic behaviour of $g$ is much better understood,
using
e.g. the result of Iwaniec\cite{Iwa2}, it is easy to show that there are at
most finitely many exceptions to conjecture 2. The difficult part of the
proof of Theorem 3 is to give an upper bound for $n$ and find properties on the
would-be-counterexample which makes it feasible to rule out these finitely many
values.

To prove our Theorem, we first note that we may choose $s=1$, since otherwise
we replace $p$ by $p'\equiv ps^{-1}\pmod{n}$. Then we have
$\frac{n}{2}+1<t<n-2$. In the first step we exclude odd values of $n$.

Assume that $n$ is an odd counterexample to Theorem 3. Define the integer $k$
by the relation
$1-\frac{1}{2^k}<\frac{t}{n}<1-\frac{1}{2^{k+1}}$, and $a := t - (1-2^{-k}) n$.
Since $n$ is odd, $2^k$ is relatively prime to $n$, hence we get
$2^k + 2^k t\bmod{n} > \frac{n}{2}$. But we have
$2^k t = (2^k-1) n + 2^k a$, hence $2^k(a+1) > \frac{n}{2}$, i.e.
$a>\frac{n}{2^{k+1}}-1$. By the definition of $k$, we have
$a<\frac{n}{2^{k+1}}$, thus $t=\left[n\left(1-\frac{1}{2^{k+1}}\right)\right]$.
Write $t=n\left(1-\frac{1}{2^{k+1}}\right)-\alpha$.

Next we give an upper bound for $2^k$. Write $t=n-b$. The cases $b=1$ and $b=2$
are excluded, since we would have $s+t=n$ resp. $2s+t=n$. If
$p\in\left[\frac{n}{2(b-1)}, \frac{n}{b}\right]$, we have
$pt\bmod n+p<\frac{n}{2}$, thus if there is some $p$ in this interval
relatively prime to $n$, we are done. Thus we have
\[
\frac{n}{b} - \frac{n}{2(b-1)} < g(n)
\]
The left hand side is decreasing with $b$, thus if $b<\sqrt{n}$ the left hand
side is at least $\frac{n(\sqrt{n}-2)}{\sqrt{n}(\sqrt{n}-1)}$, and for
$n>10000$ this is $> \frac{\sqrt{n}}{3}$. Hence we obtain the bound
$\sqrt{n}<3g(n)$. By Lemma 4 this implies
$\omega(n)\leq 4$, thus $g(n)\leq 10$ and $n<300$. Thus we may suppose
$b>\sqrt{n}$.

Let $q<2^{k+1}$ be an odd prime, and define the integer $l$ by the relation
$2^l<q<2^{l+1}$. Assume that $q\not|n$. Then $(q2^{k-l}, n)=1$, thus we get
$q2^{k-l} t\bmod{n} + q2^{k-l} > \frac{n}{2}$. Using the relation
$t=n\left(1-\frac{1}{2^{k+1}}\right)-\alpha$ with $0<\alpha<1$, this becomes
\begin{eqnarray*}
q2^{k-l} t\bmod{n} + q2^{k-l} & > & \frac{n}{2}\\
n-\frac{q n}{2^{l+1}} - q2^{k-l}\alpha + q2^{k-l} & > & \frac{n}{2}\\
\frac{n}{2}-\frac{q n}{2^{l+1}} + q2^{k-l} & > & 0
\end{eqnarray*}
Since $q\geq 2^l+1$, this implies
\[
0< -\frac{n}{2^{l+1}} + q2^{k-l}\leq -\frac{n}{2^{l+1}} + 2^{k+1}\leq
-\frac{n}{2^{l+1}} + \sqrt{n}
\]
hence $2^{l+1}\geq\sqrt{n}$. Thus $n$ is divisible by all odd primes
$\leq\sqrt{n}$. Using the elementary bound $\theta(n)>n/2$, where
$\theta(x)=\sum_{p\leq x}\log p$, this implies
$2n> e^{\sqrt{n}/2}$, which in turn implies $n<121$. However, Theorem 3 is true
for all $n<10000$, thus we conclude that it is true for all odd $n$.

Thus assume that $(n, t)$ is a counterexample to Theorem 3 with $n$ even.

We show that $t$ cannot be too close to $n/2$ or to $n$. The
proofs for these two cases run parallel, and we will only give the first one.
Set $t=\frac{n}{2} + b$. Let $p$ be any integer relatively prime to $n$,
in particular, $p$ is odd. Then we have
\[
pt = \frac{pn}{2} + bp \equiv -\frac{n}{2} + bp\pmod{n}
\]
thus if $n$ is a counterexample to our Theorem, we conclude that
$bp\not\in[n/2, 3n/2 - p]$, i.e. $p\not\in\left[\frac{n}{2b},
\frac{3n}{2b}-\frac{p}{b}\right]$. The case $b=1$ is excluded, thus the upper
bound of this interval is $\geq\frac{n}{b}$, thus in particular we have
$p\not\in\left[\frac{n}{2b}, \frac{n}{b}\right]$. But the only conditions
imposed on $p$ were that $p$ is odd and coprime to $n$. Since all even
integers are not coprime to $n$, we
get that the interval $\left[\frac{n}{2b}, \frac{n}{b}\right]$ contains no
integer relatively prime to $n$. Hence $g(n)>\frac{n}{2b}$, thus
$b>\frac{n}{2g(n)}$, i.e. $t>n/2+\frac{n}{2g(n)}$. In the same way we have
$t<n-\frac{n}{2g(n)}$.

Set $w=(t, n)$. As
$p$ runs over all integers relatively prime to $n$, $pt$ runs over all integers
with $(pt, n)=w$, and $pt\bmod{n}$ has period $n/w$. Hence there is some
$p<n/w$, relatively prime to $n$ with $pt\equiv w\pmod{n}$. But then
$pt\bmod{n}+p \leq w + n/w$, and this is $\leq n/2$, unless $w=1, 2, n/2$
or $n$. The last two cases are trivially excluded.  Thus we are left with the
cases $w=1, 2$. Now $\frac{t}{n}$
is a rational number with denominator $>\sqrt{n}$, thus applying Dirichlet's
Theorem we find an integer $d\leq\sqrt{n}$ and some $e\leq d$, such that
$\left|\frac{dt}{n} - e\right|<\frac{1}{\sqrt{n}}$.

Assume that $d=1$. Then $\left|\frac{t}{n} - e\right|<\frac{1}{\sqrt{n}}$, and
because $n/2<t<n$, we conclude $t>n-\sqrt{n}$. Together with the bound proved
above we obtain the inequality $\sqrt{n}>\frac{n}{2g(n)}$, i.e.
$2g(n)>\sqrt{n}$. Using the first statement of Lemma 4, this yields
$\omega(n)\leq 4$, thus $n<1156$, but for $n<10000$ the Theorem is already
proven. In the same way we exclude the case $d=2$.  Now
assume $d>2$. Then by Lemma 4, statement 3, we find some $a$ relatively
prime to $d$ with $\frac{d}{12}<a<\frac{5d}{12}$. Let $p$ be an integer
relatively prime to $n$ which also satisfies $p\equiv ae^{-1}\pmod{d}$. Note
that the right hand
side exists, since $(e, d)=1$. Write $p=k d + a'$. Then we have
\[
pt = \frac{pen}{d} + \theta \frac{p \sqrt{n}}{d} =  ken +\frac{a'en}{d} +
\theta\frac{p\sqrt{n}}{d}\equiv\frac{an}{d} + \theta\frac{p\sqrt{n}}{d}\pmod{n}
\]
where $\theta$ is some real number of absolute value $<1$. But $pt\bmod n$ is
$>\frac{n}{2}-p$, thus either the right hand side is $>\frac{n}{2}-p$, which
yields
\[
\frac{an}{d} + \frac{p\sqrt{n}}{d} > \frac{n}{2} - p
\]
or the right hand side is negative, which yields
\[
\frac{an}{d} - \frac{p\sqrt{n}}{d} < 0
\]
From now on, we will only consider the first inequality, because the second one
can be dealt with similarly, but gives a little stronger bounds. By the choice
of $a$ we have $a/d\leq 5/12$, thus we get $p(\frac{\sqrt{n}}{d}+1) > n/12$. By
Lemma 4, statement 2, $p$ can be chosen to be $\leq d(g(n)+1)$. Thus we obtain
the inequality $(\sqrt{n}+d)(g(n)+1)> n/12$. Since $d\leq\sqrt{n}$, we finally
conclude $g(n)>\sqrt{n}/24 - 1$. The bound $g(n)<2^{\omega(n)}$ shows that this
is only possible for $\omega(n)\leq 9$. Now the improved bound
$g(n)\leq \omega(n)^2$ lowers the bound to 7, and we can use the fifth
statement from Lemma 4 to conclude $n<(24\cdot 27)^2$, thus $\omega(n)\leq 6$
and $n<(24\cdot 23)^2=304704$. 

Assume that $p$ is some prime number, such that the least positive residue of
$ep\pmod{d}$ is in the interval $[d/12, 5d/12]$. Then by the argument above,
we get $p(\frac{\sqrt{n}}{d}+1) > n/12$ or $p|n$. Hence all primes $p$ which
satisfy this congruence condition, have to divide $n$. By the bounds given
above, it suffices to find 7 such primes to exclude the pair $(n, d)$.

To finish the proof of Theorem 3, note first that $d\leq\sqrt{304704} = 552$.
Choose some $d$, and compute $p_{\max} = \frac{10000}{100/d+1}$. Count the
number of residue classes $a$ relatively prime to $d$, with $d/12<a<5d/12$, and
call this number $N$.Count the prime numbers up to $p_{\max}$ in all reduced
residue classes $\pmod{d}$, and choose those $N$ sequences with the least
number of primes in it. If $n$ is a counterexample to Theorem 3, and $d$ is
corresponding in the sense described above, then $n$ is divisible by all these
prime numbers, in particular there are at most 6 such primes. 

Doing this for all $d\leq 552$, we found no $d$ such that there could
correspond some $n$ giving a counterexample to Theorem 3.

All computations were performed on a Silicon Graphics Indy workstation using
Mathematica 3.0.

\end{document}